\documentstyle{article}
\newcommand{\cald}{{\cal D}}
\newcommand{\calb}{{\cal B}}
\newcommand{\tpi}{T_{\pi}}
\newcommand{\xj}{x_{J}^{2}}
\newcommand{\tppi}{T_{p, \pi}}
\newcommand{\pivon}[1]{\pi(#1)}
\newcommand{\sppi}{S_{p, \pi}}
\newcommand{\lip}{\Lambda_{\left( \frac{1}{p}-1 \right)}}
\newcommand{\hochx}{^{1-\frac{1}{p}}}
\newcommand{\hochy}{^{2\left( \frac{1}{p}- \frac{1}{2}\right)}}
\newcommand{\cac}[1]{{\rm CC}(#1)}
\newcommand{\permut}{\pi : \cald \rightarrow \cald}
\newcommand{\cacp}[1]{{\rm CC_{p}}(#1)}
\newcommand{\ausdr}{\sum_{L \in max \pi (D(I))} |L|\hochy }
\newtheorem{theorem}{Theorem}
\newtheorem{lemma}{Lemma}
\newtheorem{claim}{Claim}
\newtheorem{definition}{Definition}
\begin{document}
\title{Permutations of the Haar system}
\author{Paul F.X. M\"{u}ller\thanks{Supported by E. Schr\"{o}dinger
auslandsstipendium PR.Nr J0458-PHY}\\Institut f\"{u}r Mathematik,
J. Kepler Universit\"{a}t\\ Linz, Austria\\and\\Department of
Theoretical Mathematics\\The Weizmann Institute of Science\\Rehovot,
Israel}
\maketitle
\begin{abstract}
General permutations acting on the Haar system are investigated.
We give a necessary and sufficient condition for permutations
 to induce
an isomorphism on dyadic BMO. Extensions of this
characterization to Lipschitz spaces $\lip, (0<p\leq1)$
are obtained.
When specialized to permutations
which act on one level of the Haar system only, our approach
leads to a short straightforward proof of a result due to
E.M.Semyonov and B.Stoeckert.
\end{abstract}

Let us briefly describe the setting in which we are working.
${\cal D}$ denotes the set of all dyadic intervals contained in the unit
interval.
$ \pi :{\cal D}\rightarrow{\cal D}$ denotes a permutation of the dyadic
intervals.
The operator induced by $\pi$ is determined by the equation
\begin{displaymath} T_{\pi}h_{I}=h_{\pi(I)}\end{displaymath} where
$ h_{I}$ denotes the $L_{\infty}$-normalised
 Haar function supported on the dyadic intervall I.
The main result of this paper treats  {\em general} permutations
on BMO and on Lipschitz spaces.
The condition on $\pi $ which controlles the boundedness of $T_{\pi}$
is given in termes of the Carleson constant of collections of
dyadic intervals. The proof of the general result
 given below is quite complicated.
We start therefore by
considering first a special class of permutation operators.
To study these operators on $ L_{p}$ E.M. Semyonov introduced the
parameter,
\begin{displaymath}
K=\sup \left \{\frac{|\pi^{-1}({\cal B})^{\ast}|}{|{\cal B}^{\ast}|}:
{\cal B}\subseteq{\cal D}\right\}
\end{displaymath}
where for example ${\cal B}^{*}$ denotes the pointset coverd
by the collection ${\cal B}$
\par\vspace{5 mm}
E.M Semyonov and B.St\"ockert proved the following result.
\begin{theorem}
If for every $I\in{\cal D}$ we have $|\pi(I)|=|I|$ and $K<\infty$   then
for $2\leq p<\infty$ the operator $T_{\pi}$ is bounded on $L_{p}$
\end{theorem}
We will obtain this result from
\begin{theorem}
If for every $I\in{\cal D}$ we have $|\pi(I)|=|I|$ and $K<\infty$ then
the operator $T_{\pi}$ is bounded on dyadic-BMO
\end{theorem}
{\bf Proof}: Recall that for formal series
 $f=\sum_{I\in{\cal D}}a_{I}h_{I}$ the dyadic -
BMO norm is given by
\begin{displaymath}
\left(\sup_{{\cal B}
}\frac{1}{|{\cal B}^{\ast}|}
\sum_{I\in{\cal B}}a_{I}^{2}|I|\right)^{\frac{1}{2}}
\end{displaymath}
where the supremum is extended over all collections of dyadic intervals
${\cal B}.$
We fix now  $x= \sum x_I h_I $
and obtain $T_{\pi}x = \sum x_{\pi^{-1}(I)}h_{I}$. Then choose
${\cal B} \subseteq {\cal D}$ such that
\begin{displaymath}
\frac{1}{2}||T_{\pi}x||_{BMO}^{2}\leq \frac{1}{|{\cal B}^{\ast}|}
                   \sum_{\cal B} x_{\pi^{-1}(I)}^{2} |I|
\end{displaymath}
By hypothesis this expression equals with
\begin{displaymath}  \frac{1}{|{\cal B}^{\ast}|}
                   \sum_{\cal B} x_{\pi^{-1}(I)}^{2} |\pi^{-1}(I)|
\end{displaymath}
Which equals trivially with
\begin{displaymath}\frac{|\pi^{-1}({\cal B})^{\ast}|}
{|{\cal B}^{\ast}|}
  \frac{1}{|\pi^{-1}({\cal B})^{\ast}|}
                   \sum_{\cal B} x_{\pi^{-1}(I)}^{2} |\pi^{-1}(I)|
\end{displaymath}
The last expression is of course bounded by $K ||x||_{BMO}^{2}$.
This finishes the proof of Theorem2.\par
{\bf Remark:} 1) For every permutation which satisfies
$|\pi(I)|=|(I)|$
there exist ${\cal E}\subseteq {\cal D}$ and $x\in BMO$ for which the
above chain of inequalities can be reversed. Hence for
such permutations
the condition $K<\infty$
{\em is implied by} the boundedness of $T_{\pi}$.
 \par 2)
As $T_{\pi}$ is bounded on $L_{2}$ we obtain from
\cite{JJ} Corollary2 p60
that $\tpi$ is bounded on $L_{p}$ for $2<p<\infty$. We
 thus obtained
Theorem 1
from Theorem 2 by interpolation.
\par  \par
Up to this point we
considered permutations which act on one level of
the Haar system only. We now turn to arbitrary permutations.
\par To do this we need a scale invariant measure for the size of
collections of dyadic Intervals $\calb$ :
the so called Carleson condition.
This notion was studied and carefully
analyzed by P.W. Jones in his
work on the uniform approximation property
of BMO. \par
\begin{definition}
\begin{enumerate}
\item
 $\calb$ is said to
satisfy the {\em K-Carleson condition } if
\begin{displaymath}
\sup_{J} \frac{1}{|I|}
\sum_{\{J\in \calb : J\subset I\} }|J|\leq K
\end{displaymath}

The infimum over all such K is called
the {\em Carleson constant of $\calb$ }
and will be denoted $\cac{\calb}$
\item $max\calb$ denotes those intervals in
$\calb$ which are not contained
in other intervals of $\calb$. As dyadic
intervals are nested, the collection
$max\calb$ contains only pairwise disjoint dyadic intervals.
\end{enumerate}
\end{definition}
\begin{theorem} Let $\permut$ be any permutation.
The operator $\tpi$ is
an isomorphism on BMO if and only if there exists
$M\in {\bf R}^{+}$
 such that for all $
\calb\subset\cald$
$$ \frac{1}{M}\cac{\calb}\leq\cac{\pivon{\calb}}
\leq \cac{\calb}M $$
\end{theorem}
We will prove Theorem 3 by
decomposing $\cald$ so as to control
the norm of $\tpi$ on smaller parts.
During the decomposition process
we will collect additional information
concerning the interaction of the
small pieces. Each iteration step is based on
\begin{lemma} Let $\permut$ be a permutation which satisfies
the condition of Theorem 3 Then for any $D(I)\subset \cald$
with $D(I)^* \subset I$
 ,$K\in {\bf R}^{+}$, and $x\in$ BMO we obtain a
decomposition $ D(I) = G(I) \cup S(I) $
\begin{enumerate}
\item where G(I) satisfies $$
\sum_{J\in G(I)}\xj|\pivon{J}| \leq|\pivon{D(I)}^*|
||x||_{BMO}^{2}K  $$
\item and S(I) satisfies $$
\sum_{J \in maxS(I)} \frac{|J|}{|I|} \leq \frac{M}{K}
      $$
\end{enumerate}
Consequently for  ${\cal N}(I) =\pi^{-1}(max \pivon{D(I)}$ and
${\cal O}(I) = ({\cal N}(I)\cap S(I))\cup maxS(I)$ we obtain
 $$ \sum_{J\in {\cal O}(I)}
\frac{|J|}{|I|} \leq \frac{M(M+1)}{K} $$
\end{lemma}
{\bf Proof:} We define$$S(I)=\left\{J\in D(I):
\frac{|\pivon{J}|}{|\pivon{D(I)}^{*}|}
\geq K\frac{|J|}{|I|}\right\}$$
and let $G(I)=D(I)\setminus S(I).$ The defining inequality for
$J \in G(I)$
implies that \begin{displaymath}
   \sum_{J \in G(I)} \xj |\pivon{J}| \leq
  \sum_{J \in G(I)} \xj \frac{|J|}{|I|}
K|\pivon{D(I)}^{*}| \end{displaymath}
This expression has $||x||_{BMO}
^{2} K |\pivon{D(I)}^{*}|$ as upper bound.
It remains to analyze S(I): by definition of S(I) we get
\begin{displaymath} \sum_{J\in maxS(I)} \frac{|J|}{|I|} \leq
\frac{1}{K} \sum_{J \in maxS(I)}
\frac{|\pivon{J}|}{|\pivon{D(I)}^{*}|}
\end{displaymath}
The sum on the right hand side is
dominated by the Carleson constant
of ${\pivon{maxS(I)}}$,
which by assumtion is bounded by M times the Carleson constant  of
 ${maxS(I)}$ . However ${maxS(I)}$ beeing a collection of pairwise
satisfies the 1 Carleson condition. Summing up we obtain
 \begin{displaymath} \sum_{J\in maxS(I)}
 \frac{|J|}{|I|}
\leq \frac{M}{K}  \end{displaymath}
The collection
$max\pivon{D(I)}$ like any other
 collection of disjoint dyadic intervals
has Carleson constant equal to one. Hence
$ \pi^{-1}(max\pivon{D(I)})$
satisfies the M Carleson condition.
In particular for $ L\in maxS(I)$
\begin{displaymath}
\sum_{\{J\in {\cal N}(I): J \subset L\}}
|J| \leq M|L|
\end{displaymath} Let us combine this
estimate with the previous
analysis of S(I) to describe the size of ${\cal O}(I)$.
\begin{displaymath} \sum_{J\in {\cal O}(I)}|J|
\leq \sum_{L \in max S(I)}
\sum_{\{J\in {\cal N}(I): J\subset L\}}
|J|
\leq
\sum_{L \in maxS(I)} M|L|  \leq \frac{M^{2}}{K}|I|
\end{displaymath}
The generations of the index set
${\cal O}(I)$ are used to form a stopping time
decomposition of S(I).

\begin{definition}
\begin{enumerate}
\item We first recall how generations are formed: let $G_{0}
({\cal O}(I))$ be $max{\cal O}(I)$.
Having defined $G_{0}({\cal O}(I))
,\ldots, G_{l}({\cal O}(I))$ we put
$$G_{l+1}({\cal O}(I))=
max\left({\cal O}(I) \setminus
\bigcup_{k\leq l} G_{k}({\cal O}(I))
\right) $$
\item  We now form the crucial decomposition of S(I):

 For $ k \in \{0,1,2,\ldots \}$ and $ L\in G_{k}({\cal O}(I))$
we define $D(L)=\{ J\in S(I):J \subset L\}
\setminus \{ J \subset P:
P \in G_{k+1}({\cal O}(I)) \}$ \end{enumerate} \end{definition}
{\bf Comment}
1) Consider the following identity
$$
 \sum_{J \in D(I)} \xj |\pi (J)| =
\sum_{J \in G(I)} \xj |\pi (J)|+ \sum_{L \in {\cal O}(I)}
\sum_{J \in D(L)} \xj |\pi(J)|
 $$
In view of Lemma1
the sum indexed by G(I) admits a good upper bound. We
shifted the bad behavior of the permutation into the sum indexed by
${\cal O}(I)$
However we used the hypothesis to show that this index set
is geometrically
small compared to I. This remark indicates that Lemma1 permits us to
show that a repeated application of the identity defines a converging
algorithm. \par 2)  One's first idea might be to choose maxS(I) as
index set ${\cal O}(I)$. However when one tries to prove convergence
for the associated decomposition procedure one meets serious technical
difficulties. The way to get around these complications is to choose
an index set which contains information about $maxS(I)$ {\bf and }
 $\pi^{-1}
(max\pi (D(I))$.  \par {\bf Proof of Theorem3}
The necessity of our condition is implied by the following
relation between the Carleson condition and BMO : $\cac{\calb}=
||\sum_{I\in \calb} h_{I}||_{BMO}
^{2} $. The rest of the paper is
used to show that the condition of Theorem3 is also sufficient.
We first choose
$ J_{0} \in \cald $  such that  $\frac{1}{2}||\tpi x|| ^{2}$
is bounded by
\begin{displaymath}
\frac{1}{|\pi (J_{0})|}
\sum_{\{J: \pi (J) \subset \pi (J_{0})\}}\xj |\pi (J)|
\end{displaymath}
We now let
$ \calb = \pi ^{-1}(\{J: \pi (J) \subset \pi (J_{0})\})$
and
${\cal O}_{0} = max\calb$
For $I\in {\cal O}_{0}$
we put $D(I) =\{J\in \calb: J\subset I\}$
having produced ${\cal O}_{1}\ldots {\cal O}_{l},\cald_{1}\ldots
\cald_{l}$
and ${\cal N}_{1}\ldots {\cal N}_{l}$ we choose $I\in {\cal O}_{l}$
and $D(I)\in \cald_{l}$. Lemma1 is now applied to D(I) and we obtain
G(I), S(I), $ {\cal N}(I)$ ,
and ${\cal O}(I)$. Finally S(I) is decomposed
according to the Definition2 . This gives us $\{D(L)\subset S(I):
L \in {\cal O}(I)\}$. Doing this for each $I\in {\cal O}_{l}$
allowes us to define
\begin{displaymath} {
\cal O}_{l+1} =\bigcup_{I\in {\cal O}_{l}}{\cal O}(I) ,
\cald_{l+1} =
\bigcup_{I \in {\cal O}_{l+1}} D(I), {\cal N}_{l+1}=
\bigcup_{I\in {\cal O}_{l}}{\cal N}(I)
\end{displaymath}
We thereby completed the induction step. The next two claims describe
the behaviour of our construction.
From now on, K denotes any number bigger than $4M^{2}$

\begin{claim}$ {\cal O}=\bigcup_{l\in {\bf N}}{\cal O}_{l}$ satisfies
the 2 Carleson condition
\end{claim}
{\bf Proof:} Fix $ I\in {\cal O}_{k_{0}}$. If $J\in {\cal O}_{k}$,
and $J \subset I$ then by construction we obtain $k \geq k_{0}$.
Applying the estimates of Lemma1  we obtain for $k\geq k_{0}$
\begin{displaymath}
\sum_{\{J\in {\cal O}_{k}:J\subset I\}}|J| \leq \left( \frac{M(M+1)}{K}
\right)^{k-k_{0}}|I|
\end{displaymath}
For $I\in {\cal O}_{k_{0}}$
\begin{displaymath}
\sum_{\{J \in {\cal O}: J \subset I\}}|J| =\sum_{k\geq k_{0}}
\sum_{\{J\in {\cal O}_{k}: J\subset I\}} |J|
\end{displaymath}
Invoking the observation above we obtain for this sum the following
majorization:
\begin{displaymath}
|I|\sum_{k\geq k_{0}}\left(\frac{M(M+1)}{K}\right)^{k-k_{0}}
\leq |I|2
\end{displaymath}
\begin{claim}${\cal N} =
\bigcup_{l\in {\bf N}}{\cal N}_{l}$ satisfies the
(3M) Carleson condition.
\end{claim}
{\bf Proof:} When we look back at the construction we see that each
dyadic interval lies in at most one of the collections ${\cal N}(J)$.
Fix now $I\in {\cal N}$.
Hence $I \in {\cal N}_{k_{0}}$ for some $k_{0}$
, and there exists exactly one dyadic interval $P$ such that
 $I\in {\cal N}(P)$. This remark gives the representation
\begin{displaymath}
\{J\in {\cal N}: J\subset I\}=\{J\in {\cal N}(P): J\subset P\}\cup
\{{\cal N}(L): L\in {\cal O}, L\subset P\}
\end{displaymath}
We thus obtain the following identity:
\begin{displaymath}
\frac{1}{|I|}\sum_{\{J\in {\cal N}: J\subset I\}}|J| =
\frac{1}{|I|}\sum_{J\in {\cal N}(P)}|J|+
\frac{1}{|I|}\sum_{\{L\in {\cal O}:L\subset P\}}
 \sum_{\{J\in {\cal N}(L)\}} |J|
\end{displaymath}
The first summand is simply estimated by $\cac{{\cal N}(P)}$
which in turn
is less than $M$.The second term is majorized by
\begin{displaymath}
\frac{1}{|I|}
\sum_{\{L\in {\cal O}:L \subset P\}} M|L| \leq M \cac{{\cal O}}
\end{displaymath}
Having proved claim2,{\bf we resume the proof of Theorem3.} \par
We know now that ${\cal N}$ satisfies the 3M Carleson condition.
Observe that
$\pivon{{\cal N}} = \bigcup_{I \in {\cal O}}max\pivon{D(I)}$
Therefore
\begin{displaymath}
\frac{1}{|\pi (J_{0})|}
\sum_{I\in {\cal O}} |max \pi(D(I))^{*}| \leq
\cac{\pi({\cal N})}
\end{displaymath}
By hypothesis on
$\pi$ we get $\cac{\pi({\cal N})}\leq M\cac{{\cal N}}$.
Recall next that
$\{G(I): I\in {\cal O}\}$ is a decomposition of $\calb$.
We thus obtain the following final estimates.
$$
\frac{1}{2} || \tpi x||_{BMO}
^{2} \leq \frac{1}{|\pi (J_{0})|}
\sum_{\{\pi(J): J\in \calb\}} \xj |\pi(J)| =
\frac{1}{|\pi(J_{0})|}\sum_{I \in {\cal O}}
\sum_{\{\pi(J): J\in G(I)\}}\xj |\pi(J)|  $$
$$ \leq \frac{1}{|\pi(J_{0})|}
\sum_{I\in {\cal O} } |max\pi(D(I)^{*}|||x||^{2}  $$
As we observed above this sum is bounded from above
by  $ 3M^{2} ||x||_{BM}^{2}$

This proves that $\tpi $ is a bounded operator on $BMO$. As
the hypothesis
of the theorem is symmetric in $\pi$ and $\pi^{-1}$ we conclude that
$T_{\pi^-1} =\tpi^{-1}$ is bounded as well.
\par {\bf Extensions to Lipschitz functions:}
The dyadic BMO condition appears as a natural limit of some
Lipschitz condition for martingales. We describe subsequenty an
extension of our main result to Lipschitz spaces. l
For $ \varphi \in L^2 $ say and
$0<p\leq 1$ the Lipschitz condition assumes the form
\begin{displaymath} ||\varphi||_{\lip}=
\sup_{I\in \cald }\left(
\int_{I} |\varphi - \varphi_{I}|^{2}\frac{dt}{|I|}\right)
^{\frac{1}{2}} {|I|\hochx}
\end{displaymath}
where $ \varphi_{I} $ denotes the meanvalue of $\varphi $ over $I$.
The particular interest in these class of functions stems from
a duality relation due to C. Herz,\cite{HE},
 which generalizes C.Feffermans
duality theorem. Herz's theorem identifies $\lip $ with the
dual space of dyadic $H^p (0<p\leq 1).$ \par
We shall now discuss permutations of Haar functions which
are normalised in $\lip$.
The norm of the $L^{\infty}$ normalised functions $h_I$ in $\lip$
equals $|I|\hochx$.
Hence the operator induced by the permutation
$\permut$ is given by  the relation
\begin{displaymath}T_{p,\pi}\left( \frac{h_I}{|I|\hochx} \right) =
\frac{h_{\pi (I)}}{|\pi (I)|\hochx}
\end{displaymath}
It is useful to observe that the $\lip$ norm of
$ f= \sum a_{I}\frac{h_{I}}{|I|\hochx} $
can be expressed in terms of the coefficients $a_{I}$. In fact
we obtain
\begin{displaymath}
||f||_{\lip}^{2} = \sup_{\calb}\left(\frac{1}{|I|\hochy}
\sum_{J \subset I} a_{J}^{2}|J|\hochy \right)
\end{displaymath}
This formula suggests how to extend properly the notion of
Carleson - constant.
\begin{definition} For a collection $\calb$ of dyadic intervals,
the {\em Carleson p-constant} is given by
\begin{displaymath}
\sup_{I \in \calb}\left( \frac{1}{|I|\hochy }\sum_{\{J\in \calb:
J \subseteq I \}} |J|\hochy \right)
\end{displaymath}
It will be denoted by $ {\rm CC_{p}}(\calb) $
\end{definition}
Our extension of Theorem 3 reads now as follows
\begin{theorem} Let $\permut$ be a permutation.  Let $0<p\leq 1 $.
The operator $T_{p,\pi} $ is an isomorphism on $\lip$ if and
only if there exists $M \in {\bf R^+} $ such that for every
$\calb \subseteq \cald $
\begin{displaymath}
\frac{1}{M} {\rm CC_p}(\calb) \leq {\rm CC_p}(\pi(\calb)) \leq
{\rm CC_p}(\calb) M
\end{displaymath}
\end{theorem}
The proof of theorem4 follows the pattern explained above.
The first step is again given by the following
\begin{lemma} Let $\permut$ be a permutation which satisfies
the condition of Theorem 4 Then for any $D(I)\subset \cald$
with $D(I)^* \subset I$
 ,$K\in {\bf R}^{+}$, and $x\in \lip$ we obtain a
decomposition $ D(I) = G(I) \cup S(I) $
\begin{enumerate}
\item where G(I) satisfies $$
\sum_{\{\pivon{J}:J\in G(I)\}}\xj|\pivon{J}|\hochy
\leq\ausdr||x||_{\lip}^{2}K
 $$
\item and S(I) satisfies $$
\sum_{J\in maxS(I)}\left( \frac{|J|}{|I|} \right)\hochy
\leq \frac{M}{K}
      $$
\end{enumerate}
Consequently for  ${\cal N}(I) =\pi^{-1}(max \pivon{D(I)}$ and
${\cal O}(I) = ({\cal N}(I)\cap S(I))\cup maxS(I)$ we obtain
 $$ \sum_{J\in {\cal O}(I)}\left(
\frac{|J|}{|I|} \right)\hochy\leq \frac{M(M+1)}{K} $$
\end{lemma}
{\bf Proof:}
We define$$S(I)=\left\{J\in D(I):\frac{|\pivon{J}|\hochy}{\ausdr}
\geq K\left(\frac{|J|}{|I|}\right)\hochy\right\}$$
and let $G(I)=D(I)\setminus S(I).$ The defining inequality for
$J \in G(I)$
implies that \begin{displaymath}
   \sum_{J \in G(I)} \xj |\pivon{J}|\hochy \leq
  \sum_{J \in G(I)} \xj\left( \frac{|J|}{|I|} \right)\hochy
\ausdr\end{displaymath}
This expression has $$||x||_{\lip}^{2} K \ausdr$$ as upper bound.
It remains to analyze S(I): by definition of S(I) we get
\begin{displaymath} \sum_{J\in maxS(I)} \left(
 \frac{|J|}{|I|}\right)\hochy \leq
\frac{1}{K} \sum_{J \in maxS(I)} \frac{|\pi(J)|
\hochy}{\ausdr}
\end{displaymath}
The sum on the right hand side is
dominated by the Carleson p constant
of ${\pivon{maxS(I)}}$,
which by assumtion is bounded by M times the Carleson p constant  of
 ${maxS(I)}$ . This in turn is bounded by one. We obtained:
 \begin{displaymath} \sum_{J\in maxS(I)}\left( \frac{|J|}{|I|}
\right)\hochy
\leq \frac{M}{K}  \end{displaymath}
The collection
$max\pivon{D(I)}$
has Carleson p constant equal to one. \par
Hence $ \pi^{-1}(max\pivon{D(I)})$
satisfies the M Carleson p
condition. In particular for $ L\in maxS(I)$
\begin{displaymath}
\sum_{\{J\in {\cal N}(I): J \subset L\}} |J|\hochy \leq M|L|\hochy
\end{displaymath} This estimate and the previous
analysis of S(I) gives an upper bound for the
 size of ${\cal O}(I)$.
\begin{displaymath} \sum_{J\in {\cal O}(I)}|J|\hochy
\leq \sum_{L \in max S(I)}\sum_{\{J\in {\cal N}(I): J\subset L\}}
|J|\hochy  \end{displaymath} \begin{displaymath}
\leq
\sum_{L \in maxS(I)} M|L|\hochy  \leq \frac{M^{2}}{K}|I|
\end{displaymath}

{\bf Proof of Theorem 4}
The necessity of our condition is implied by the following
fact $$\cacp{\calb}=
||\sum_{I\in \calb}\frac{ h_{I}}{|I|\hochx}||_{\lip}^{2}. $$
We show now that the condition of Theorem4 is also sufficient.
We first choose
$ J_{0} \in \cald $
such that  $\frac{1}{2}||\tppi x||_{\lip} ^{2}$
is bounded by
\begin{displaymath}
\frac{1}{|\pi (J_{0})|\hochy}
\sum_{\{J: \pi (J) \subset \pi (J_{0})\}}\xj |\pi (J)|\hochy
\end{displaymath}
We now let $ \calb = \pi ^{-1}(\{J: \pi (J) \subset \pi (J_{0})\})$.
Using Lemma2 we form as in the proof of Theorem3 the collections
${\cal O},{\cal N} $, $\{G(I): I\in {\cal O}\}$ and $\{D(I): I \in
{\cal O}\}$  such that :
\begin{enumerate}
\item $\cacp{{\cal N}} \leq 3M$
\item $\{G(I): I \in {\cal O}\}$ is a decomposition of $\calb$
\item $$ \sum_{J \in G(I)} \xj |\pi(J)|\hochy \leq
\ausdr ||x||_{\lip}^{2}K$$
\end{enumerate}
As $\pivon{{\cal N}} = \bigcup_{I \in {\cal O}}max\pivon{D(I)}$
we obtain
\begin{displaymath}
\frac{1}{|\pi (J_{0})|\hochy}
\sum_{I\in {\cal O}} \ausdr \leq
\cacp{\pi({\cal N})}
\end{displaymath}
By hypothesis on $\pi$ we get $\cacp{\pi({\cal N})}\leq
M\cacp{{\cal N}}$.
Moreover
$\{G(I): I\in {\cal O}\}$ is a decomposition of $\calb$.
We thus obtain
$$
\frac{1}{2} || \tppi x||_{\lip}^{2} \leq \frac{1}
{|\pi (J_{0})|\hochy}
\sum_{\{\pi(J): J\in \calb\}} \xj |\pi(J)|\hochy = $$
$$ \frac{1}{|\pi(J_{0})|\hochy}\sum_{I \in {\cal O}}
\sum_{\{\pi(J): J\in G(I)\}}\xj |\pi(J)|  $$
$$ \leq \frac{1}{|\pi(J_{0})|\hochy}
\sum_{I\in {\cal O} } \ausdr||x||_{\lip}^{2}  $$
We observed already that this sum is bounded
by  $ 3M^{2} ||x||_{\lip}^{2}$

This proves that $\tppi $ is a bounded operator on $\lip$.
Repeating this procedure with $\pi$ replaced by $\pi^{-1}$
wwe are able to conclude that
$T_{p,\pi^-1} =\tppi^{-1}$ is bounded as well.
\par {\bf Remark}  \par
In view of the above mentioned duality results due
to C. Herz and C. Fefferman the transposed
operator of $\tppi$, $(0<p\leq 1)$ is given
by the operator $\sppi$ defined by $\pi^{-1}$ acting
on $H^p$ normalised Haar functions. More precisely
$ \sppi :H^p \rightarrow H^p $ is given by the relation.
$$ \sppi \left( \frac{h_{I}}{|I|^{\frac{1}{p}}}\right) =
\frac{h_{\pi^{-1}(I)}}{|\pi^{-1}(I)^{\frac{1}{p}}} $$
An equivalent formulation of Theorem 4 is thus given by
\begin{theorem} Let $\permut$ be a permutation.  Let $0<p\leq 1 $.
The operator $\sppi $ is an isomorphism on dyadic $H^p$ if and
only if there exists $M \in {\bf R^+} $ such that for every
$\calb \subseteq \cald $
\begin{displaymath}
\frac{1}{M} {\rm CC_p}(\calb) \leq {\rm CC_p}(\pi(\calb)) \leq
{\rm CC_p}(\calb) M
\end{displaymath}
\end{theorem}

\end{document}